\documentstyle[12pt]{article}
\newcommand{\bt}{\begin{Theorem}}
\newcommand{\et}{\end{Theorem}}
\newcommand{\bi}{\begin{itemize}}
\newcommand{\ei}{\end{itemize}}
\newcommand{\bea}{\begin{eqnarray}}
\newcommand{\eea}{\end{eqnarray}}
\newtheorem{Definition}{Definition}[section]
\newtheorem{Theorem}[Definition]{Theorem}
\newtheorem{Lemma}[Definition]{Lemma}

\newtheorem{Proposition}[Definition]{Proposition}
\newtheorem{Corollary}[Definition]{Corollary}

\newcommand{\be}{\begin{equation}}
\newcommand{\ee}{\end{equation}}

\newcommand{\R}{\mbox{$I\!\!R$}}%
\newcommand{\C}{\mbox{\rule{.1mm}{2.5mm}$\!\!C$}}%

\begin{document}

\begin{center}
{\protect\large\bf Analogues of the  Wiener-Tauberian and 
Schwartz theorems for radial functions on symmetric spaces}
\end{center}
\vskip.25in
\begin{center}
{\protect \bf E. K. Narayanan}\footnote{
The first author was supported in part by a grant 
from UGC via DSA-SAP} and {\protect \bf A. Sitaram}\footnote{
The second author was supported by IISc Mathematics Initiative}
\end{center}
\begin{abstract}
We prove a Wiener-Tauberian theorem for the $L^1$ spherical 
functions 
on a semisimple Lie group of arbitrary real rank. We also
establish a Schwartz type theorem for complex groups. As a
corollary we obtain a Wiener-Tauberian type result for 
compactly supported distributions.
\vskip.15in
{\bf AMS Classification 2000:}~Primary 43A20, 43A90, 
Secondary 43A80, 43A30.
\end{abstract}
\begin{center}
\section{\protect\large\bf Introduction}
\end{center}
\setcounter{equation}{0}

Two celebrated theorems from classical analysis dealing with
translation invariant subspaces are the Wiener-Tauberian theorem
and the Schwartz theorem.
Let $f \in L^1(\R)$ and $\tilde{f}$ be its Fourier transform. 
Then the 
celebrated Wiener-Tauberian theorem says that the ideal 
generated by $f$
is dense in $L^1(\R)$ if and only if $\tilde{f}$ is a nowhere 
vanishing 
function
on the real line. 

The result due to L. Schwartz says that, every closed
translation invariant subspace $V$ of $C^\infty(\R)$ is
generated by the exponential polynomials in $V.$ In particular,
such a $V$ contains the function $x \rightarrow e^{i\lambda x}$
for some $\lambda \in ~\C.$ Interestingly, this result fails
for $\R^n,$ if $n \geq 2.$ Even though an exact analogue of
the Schwartz theorem fails for $\R^n ~n \geq 2,$ it
follows from 
the well known theorem of Brown-Schreiber-Taylor \cite{BST} 
that, if
$V \subset C^\infty (\R^n)$ is a closed subspace which is
translation and rotation invariant then $V$ contains a $\psi_{
s}$ for some $s \in ~\C$ where $$\psi_s(x)
= \frac{J_{\frac{n}{2}-1}(s |x|)}{ (s |x|)^{
\frac{n}{2}-1}} = \int_{S^{n-1}}~e^{is x.w}~d\sigma(w).$$
Here $J_{\frac{n}{2}-1}$ is the Bessel function of the first 
kind and order $n/2 -1$ and $\sigma$ is the unique, normalized
rotation invariant measure on the sphere $S^{n-1}.$ It also
follows from the work in \cite{BST} that $V$ contains all the
exponentials $e^{z.x},$ if $z = (z_1, z_2, \dots z_n) \in ~\C^n$
satisfies $z_1^2 + z_2^2 + \cdots + z_n^2 = s^2.$

Our aim in this paper is to prove analogues of these results
in the context of non compact semisimple Lie groups.

\noindent
{\large\bf Notation and preliminaries:}~
For any unexplained terminology we refer to \cite{H}.
Let $G$ be a connected non compact 
semisimple Lie group with finite center and $K$ a fixed maximal
 compact
subgroup of $G.$ Fix an Iwasawa decomposition $G = KAN$ and 
let $\bf{a}$
be the Lie algebra of $A.$ Let $\bf {a^*}$ be the real dual 
of $\bf{a}$ and
$\bf{a^*_{~\C}}$ its complexification. Let $\rho$ be the half 
sum of 
positive
roots for the adjoint action of $\bf{a}$ on $\bf{g},$ the Lie 
algebra of 
$G.$ 
The Killing form induces a positive definite form $<. , >$
on $\bf{a^*} \times 
\bf{a^*}.$ Extend this form to a bilinear form on 
$\bf{a^*_{~\C}}.$ We will
use the same notation for the extension as well. Let $W$ be the
 Weyl 
group of 
the symmetric space $G/K.$ Then there is a natural action of 
$W$ on 
$\bf{a}, \bf{a^*}, \bf{a^*_{~\C}}$ and $<.,.>$ is invariant 
under this
action.

For each $\lambda \in \bf{a^*_{~\C}},$ let $\varphi_\lambda$ 
be the 
elementary
spherical function associated with $\lambda.$ Recall that 
$\varphi_\lambda$ is
given by the formula $$ \varphi_\lambda(x) = 
\int_{K}~e^{(i\lambda-\rho)
(H(xk))}~dk \quad x \in G .$$ 
See \cite{H} for more details. It is known that $\varphi_\lambda
= \varphi_{\lambda^{'}}$ if and only if $\lambda^{'} = 
\tau \lambda$
for some $\tau \in W.$ Let $l$ be the dimension of $\bf{a}$ and 
$F$ 
denote the set (in $~\C^l$ ) $$F = {\bf{a^*}} + i C_\rho
\quad \mbox{where}~ C_\rho = ~\mbox{convex hull of} 
~\{ s\rho: s \in W \}.$$
Then it is a well known theorem of Helgason and Johnson that 
$\varphi_\lambda$
is bounded if and only if $\lambda \in F.$ 

Let $I(G)$ be the set of all complex valued spherical functions 
on $G$,
that is $$I(G) = \{ f:~f(k_1 x k_2) = f(x): k_1, k_2 \in K, x 
\in G \}.$$
Fix a Haar measure $dx$ on $G$ and let $I_1(G) = I(G) 
\cap L^1(G).$ Then
it is well known that $I_1(G)$ is a commutative Banach algebra
 under 
convolution and that the maximal ideal space of $I_1(G)$ can be
 identified
with $F/W.$ 

For $f \in I_1(G),$ define its spherical Fourier transform, 
$\hat{f}$ on
$F$ by $$\hat{f}(\lambda) = \int_G~f(x)~\varphi_{-\lambda}(x)
~dx.$$ 
Then $\hat{f}$ is a $W$ invariant bounded function on $F$ 
which is 
holomorphic in the interior $F^0$ of $F,$ and continuous on $F.$
 Also 
$\widehat{f*g} = \hat{f} \hat{g}$ where the convolution of $f$
 and $g$ is 
defined by $$ f*g(x) = \int_G~f(xy^{-1})~g(y)~dy.$$

Next, we define the $L^1$- Schwartz space of $K$-biinvariant 
functions on
$G$ which will be denoted by $S(G).$ Let $x \in G.$ Then 
$x = k~expX,
~k \in K,~X \in \bf{p},$ where $\bf{g} = \bf{k} + \bf{p} $ 
is the Cartan
decomposition of the Lie algebra $\bf{g}$ of $G.$
Put $\sigma(x) = \|X\|,$ where
$\|.\|$ is the norm on $\bf{p}$ induced by the Killing form. 
For any left 
invariant differential operator $D$ on $G$ and any integer 
$r \geq 0,$
we define for a smooth $K$-biinvariant function $f$ 
$$p_{D, r}(f) = 
\sup_{x \in G}~(1+\sigma(x))^r~|\varphi_0(x)|^{-2}~|Df(x)|$$ 
where 
$\varphi_0$
is the elementary spherical function corresponding to 
$\lambda = 0.$ 
Define $$S(G) = \{f:~p_{D, r}(f) < \infty ~\mbox{for all}~D, 
r
\}.$$ Then $S(G)$ becomes
a Fr$\acute e$chet space when equipped with the topology induced by 
the family of
semi norms $p_{D, r}.$ 

Let $P = P(\bf{a^*_{~\C}})$ be the symmetric algebra over $
\bf{a^*_{~\C}}.$
Then each $u \in P$ gives rise to a differential operator 
$\partial(u)$ on
$\bf{a^*_{~\C}}.$
Let $Z(F)$ be the space of functions $f$ on $F$ satisfying the
following conditions: \

\indent
(i) $f$ is holomorphic in $F^0$ (interior of $F$) and 
continuous on $F,$ \

\indent
(ii) If $u \in P$ and  $m \geq 0$ is any integer, then 
$$q_{u, m}(f) = 
\sup_{\lambda \in F^0}~(1+\| \lambda \|^2)^m~|\partial(u)
f(\lambda) |
< \infty,$$

\indent
(iii) $f$ is W invariant .

Then $Z(F)$ is an algebra under pointwise multiplication and 
a Fr$\acute e$chet
space when equipped with the topology induced by the seminorms 
$q_{u, m}.$

If $a \in Z(F)$ we define the ``wave packet'' $\psi_a$ on 
$G$ by 
$$\psi_a(x) = \frac{1}{|W|}~\int_{\bf{a^*}}~a(\lambda)~
\varphi_{\lambda}(x)
~|c(\lambda)|^{-2}~d\lambda,$$ where $c(\lambda)$ is the 
well known
Harish-Chandra $c$-function. By the Plancherel theorem due to 
Harish-Chandra
we also know that the map $f \rightarrow \hat{f}$ extends to
 a unitary map
from $L^2(K \backslash G / K)$ onto $L^2({\bf{a^*}}, 
|c(\lambda)|^{-2}d\lambda).$
We are now in a position to state a result of
 Trombi-Varadarajan \cite{TV}.

\bt\label{1.1}
(i) If $f \in S(G)$ then $\hat{f} \in Z(F).$ \

\noindent
(ii) If $a \in Z(F)$ then the integral defining the
 ``wave packet'' 
$\psi_a$
converges absolutely and $\psi_a \in S(G).$ Moreover, 
$\hat{\psi_a}
= a.$ 

\noindent
(iii) The map $f \rightarrow \hat{f}$ is a topological 
linear isomorphism 
of
$S(g)$ onto $Z(F).$

\et

The plan of this paper is as follows: in the next section we
prove a Wiener-Tauberian theorem for $L^1(K \backslash G/K)$
assuming more symmetry on the generating family of functions.
In the final section we establish a Schwartz type theorem 
for complex semisimple Lie groups. As a corollary we also
obtain a Wiener-Tauberian type theorem for compactly supported
distributions on $G/K.$

\begin{center}
\section{\bf A Wiener-Tauberian theorem for $L^1(K \backslash
G / K)$}
\end{center}

In \cite{EM}, Ehrenpreis and Mautner observed 
that 
an exact analogue of the Wiener-Tauberian theorem
is not true for the commutative algebra of $K$-biinvariant 
functions on the
semisimple Lie group $SL(2, \R).$ Here  $K$ is the maximal 
compact subgroup
$SO(2).$ However, in the same paper it was also proved that 
an additional
``not too rapidly decreasing condition'' on the spherical 
Fourier transform
of a function suffices to prove an analogue of the
 Wiener-Tauberian
theorem. That is, if $f$ is a $K$-biinvariant integrable 
function on $G = 
SL(2, \R)$ and its spherical Fourier transform $\hat{f}$ 
does not vanish
anywhere on the maximal ideal space (which can be identified 
with a certain
strip on the complex plane) then the function $f$ generates 
a dense 
subalgebra of $L^1(K \backslash G / K)$ provided $\hat{f}$ 
does not 
vanish too fast at $\infty.$ See \cite{EM} for precise 
statements.

There have been a number of attempts to generalize these 
results to
$L^1(K \backslash G /K)$ or $L^1(G/K)$ where $G$ is a non 
compact 
connected semisimple Lie group with finite center. Almost
 complete results
have been obtained when $G$ is a real rank one group. We 
refer the reader
to \cite{BW}, \cite{BBHW} 
\cite {RS98} and \cite{S88} for results on rank one case. 
See also \cite{RS97}for a result on the whole group $SL(2,\R).$

In \cite{S80}, it is proved that
under suitable conditions on the spherical Fourier transform 
of a single
function $f$ an analogue of the Wiener-Tauberian theorem holds
 for $L^1(K
\backslash G/ K)$, with no assumptions on the rank of $G.$ 
Recently, the first named 
author improved this result to include the case of a family 
of functions rather than a single function (see\cite{N}). 
One difference
between rank one results and higher rank results 
has been the precise
form of the ``not too rapid decay condition''. In \cite{S80}
 and \cite{N} this
condition on the spherical Fourier transform of a function 
is assumed to be true on the whole maximal domain, while 
for rank one
groups it suffices to have this condition on ${\bf a^*}$
(see \cite{BW} and \cite{RS98} (An important corollary of
this is that, in the rank one case one can get a Wiener-
Tauberian type theorem for a wide class of functions purely
in terms of the non vanishing of the spherical Fourier transform
in a certain domain without having to check any decay conditions,
see \cite{MRSS}, Theorem 5.5).
In the first part of this paper we show that such a stronger
result is true for higher rank case as well provided we
assume more symmetry on the generating family of functions,
and again as a corollary we get a result of the type alluded
to in the parenthesis above.

If $dim {\bf a^*} = l,$ then ${\bf a}_{~\C}^*$ may be 
identified with $~\C^l$ and a point $\lambda \in 
{\bf a}_{~\C}^*$
will be denoted $\lambda = (
\lambda_1, \lambda_2, \dots \lambda_l).$ Let $B_R$
denote the ball of radius $R$ centered at the origin in
${\bf {a}^*}$ and $F_R$ denote the domain in ${\bf {a}}_{~\C}^*$
defined by $$F_R = \{ \lambda \in {\bf a}_{~\C}^*:~
\|Im (\lambda)\| < R \}.$$ 

For $a > 0,$ let $I_a$ denote the strip in the complex 
plane defined by $$I_a = \{ z \in ~\C:~|Im z| < a \}.$$ 
Now, suppose that $f$ is a holomorphic
function on $F_R$ and $f$ depends only on $(\lambda_1^2 +
\lambda_2^2 + \cdots +\lambda_l^2)^{\frac{1}{2}}.$ Then it is
easy to see that $$g(s) = f(\lambda_1, 
\lambda_2, \dots \lambda_l)$$
where $s^2 = \lambda_1^2 + \lambda_2^2 + \cdots + \lambda_l^2$
defines an even holomorphic function on $I_R$ and vice versa.

We will need the following lemmas. Let $A(I_a)$ denote the 
collection of functions $g$ with the properties:

\noindent
(i) $g$ is even, bounded and holomorphic on $I_a,$

\noindent
(ii) $g$ is continuous on $\bar{I}_a,$

\noindent
(iii) $\lim_{|s| \rightarrow \infty} g(s) = 0.$

Then $A(I_a)$ with the supremum norm 
is a Banach algebra under pointwise multiplication.

\begin{Lemma}
Let $\{ g_\alpha:~\alpha \in I\}$ be a collection of
functions in $A(I_a).$ Assume that there exists
no $s \in \bar{I}_a$ such that $g_\alpha (s) = 0~\forall
\alpha \in I.$ Further assume that there exists $\alpha_0 \in I$
such that $g_{\alpha_0}$ does not decay very rapidly on
$\R$, i.e,
$$ \limsup_{|s| \rightarrow \infty}~ |g_{\alpha_0}(s)|~
e^{ke^{|s|}} > 0$$ on $\R$ for all $k > 0.$ Then the 
closed ideal generated by $\{ g_\alpha:~\alpha \in I \}$ 
is whole of $A(I_a).$
\end{Lemma}

\noindent
{\bf Proof:}~Let $\psi$ be a suitable biholomorphic map which
maps the strip $I_a$ onto the unit disc (see \cite{BW}). Let
$h_\alpha(z) = g_\alpha(\psi(z)).$ Then $h_\alpha \in
A_0(D),$ where $A_0(D)$ is the collection of even holomorphic
functions $h$ on the unit disc, continuous up to the boundary
and $h(i) = h(-i) = 0.$ The not too rapid decay condition 
on $\R$ is precisely what is needed to apply the Beurling-Rudin 
theorem to complete the proof. We refer to \cite{BW} (see 
the proof of Theorem 1.1 and Lemma 1.2) for the details.

Let $p_t$ denote the 
$K$-biinvariant function defined by $\hat{p}_t(\lambda)
= e^{-t \langle \lambda, \lambda \rangle}.$ It is easy to
see that $p_t \in S(G).$

\begin{Lemma} Let $J \subset L^1(K \backslash G/K)$ be a closed
ideal. If $p_t \in J$ for some $t > 0,$ then $J = L^1(K 
\backslash G/K).$
\end{Lemma}

\noindent
{\bf Proof:}~ This follows from the main result in \cite{N}
or \cite{S80}.

Before we state our main theorem we define the following: We
say that a function $f \in L^1(K \backslash G /K)$ is 
{\it radial} if the spherical Fourier transform 
$\hat{f}(\lambda)$
is a function of $(\lambda_1^2 + \lambda_2^2 + \cdots 
\lambda_l^2)^{\frac{1}{2}}.$ Notice that, if the group $G$
is of real rank one, then the class of {\it radial} functions
is precisely the class of $K$-biinvariant functions in $L^1(G).$
When the group $G$ is complex, it is possible to describe the
class of {\it radial} functions (see next section). The 
following is our main theorem in this section:

\bt
Let $\{ f_\alpha:~\alpha \in I \}$ be a collection of 
{\it radial} functions in $L^1(K \backslash G /K).$ Assume
that the spherical transform $\hat{f}_\alpha$ extends as a
bounded holomorphic function to the bigger domain $F_R,$ 
where $R > \| \rho \|$ with
$\lim_{|\lambda| \rightarrow \infty}~\hat{f}_{\alpha}(\lambda)
= 0$ for all
$\alpha$ and that there exists no $\lambda \in
F_R$ such that $\hat{f}_{\alpha} (\lambda) = 0$
for all $\alpha.$ Further assume that there exists an $\alpha_0$
such that $\hat{f}_{\alpha_0}$ does not decay too rapidly on
${\bf a^*},$ i.e, 
$$\limsup_{|\lambda| \rightarrow \infty} |\hat{f}_{
\alpha_0}(\lambda)|~exp(k e^{|\lambda|}) > 0 $$ for
all $k > 0$ on ${\bf a^*}.$
Then the closed ideal generated by $\{f_\alpha:~\alpha \in I\}$
is all of $L^1(K \backslash G /K).$
\et

\noindent
{\bf Proof:}~ Since $f_\alpha$ is radial, 
each $\hat{f}_{\alpha}$ gives rise to an even 
bounded holomorphic function $g_\alpha(s)$ on the strip
$I_R.$ If $|\rho| < a < R,$ then the collection $\{
g_\alpha(s), \alpha \in I \}$ satisfies the hypotheses in
Lemma 2.1 on the domain $I_a.$ 
It follows that the family $\{g_\alpha\}$ generates $A(I_a).$
In particular, we have a sequence
$$h_1^n(s) g_{\alpha_1(n)}(s)
+ h_2^n(s) g_{\alpha_2(n)}(s) + \cdots + h_k^n(s) g_{\alpha_k
(n)}(s) \rightarrow e^{-\frac{s^2}{2}}$$ uniformly on 
 $\bar{I}_a,$ 
where $g_{\alpha_j (n
)} $ are in the given family and $h_j^n(s) \in A(I_a).$

Notice that each $h_j^n$ can be
viewed as a holomorphic function on the domain $F_a$ 
contained in ${\bf a}_{~\C}^*$ which
depends only on $(\lambda_1^2 + \lambda_2^2 + \cdots +
\lambda_l^2)^{\frac{1}{2}}.$ Since $h_j^n$ are bounded
and $|\rho| < a$ it can be easily verified that $e^{
-\frac{\langle \lambda, \lambda \rangle}{2}} h_j^n(
\lambda) \in Z(F).$
Again, an application of the
Cauchy integral formula says that 
$$e^{-\frac{\langle \lambda,
 \lambda \rangle}{2}} h_1^n(\lambda) \hat{f}_{\alpha_1(n)}
(\lambda) + e^{-\frac{\langle \lambda,
 \lambda \rangle}{2}} h_2^n(\lambda) \hat{f}_{\alpha_(n)}
(\lambda) +  \cdots
e^{-\frac{\langle \lambda,
 \lambda \rangle}{2}} h_k^n(\lambda) \hat{f}_{\alpha_k(n)}
(\lambda) $$ converges to $e^{- \langle \lambda, \lambda \rangle
}$ in the topology of $Z(F)$ (see the proof of Theorem 1.1 in 
\cite{BW}).
By Theorem 1.1 this simply means 
that the ideal generated by 
$\{ f_\alpha:~\alpha \in I \}$ in $L^1(
K \backslash G/K)$ contains the function $p$ where $\hat{p}
(\lambda) = e^{- \langle \lambda, \lambda \rangle}.$ We
finish the proof by appealing to Lemma 2.2.

\begin{Corollary}
Let $\{ f_\alpha: \alpha \in I \}$ be a family of radial
functions satisfying the hypotheses in Theorem 2.3. Then the
closed subspace spanned by the left $G-$translates of the
above family is all of $L^1(G/K).$ 
\end{Corollary}

\noindent{\bf Proof:}~Let $J$ be the closed subspace generated
by the left translates of the given family. By Theorem 2.3,
$L^1(K \backslash G/K) \subset J.$ Now, it is easy to see that
$J$ has to be equal to $L^1(G/K).$

\begin{Corollary}
Let $\{ f_\alpha: \alpha \in I \}$ be a family of $L^1-$radial
functions. Assume that each $\hat{f}_\alpha$ extends to a 
bounded
holomorphic function to the bigger domain $F_R$ for some $R
> \| \rho \|.$ Assume further that 
$\lim_{\| \lambda \| \rightarrow \infty}~
\hat{f}_\alpha (\lambda) \rightarrow 0.$ If there exists an 
$\alpha_0$
such that $f_{\alpha_0}$ is {\bf not} equal to a real analytic 
function almost everywhere, then the left $G-$translates of 
the above family span a dense subset of $L^1(G/K).$ 
\end{Corollary}

\noindent
{\bf Proof:}~This follows exactly as in Theorem 5.5 of 
\cite{MRSS}. 

\begin{center}
\section{\bf Schwartz theorem for complex groups}
\end{center}

When $G$ is a connected non compact semisimple Lie group
of real rank one with
finite center, a Schwartz type theorem was proved by Bagchi and
Sitaram in \cite{BS79}. Let $K$ be a maximal compact 
subgroup of $G,$
then the result in \cite{BS79} states the following: Let 
$V$ be a closed
subspace of $C^\infty (K \backslash G/K)$ with the property 
that $f \in V$ implies $w*f \in V$ for every compactly supported
$K$-biinvariant distribution $w$ on $G/K,$ then $V$ contains an
elementary spherical function $\varphi_\lambda$ for some $
\lambda \in {\bf a}_{~\C}^*.$ This was done by establishing
a one-one correspondence between ideals in $C^\infty(
K \backslash G /K)$ and that of $C^\infty(\R)_{even}.$ This also
proves that a similar result can not hold for higher rank 
groups.

Going back to $\R^n,$ we notice that 
if $f \in C^\infty(\R^n)$ is radial, then the translation 
invariant subspace $V_f$ generated by $f$ is also rotation 
invariant. It follows from \cite{BST} that $V_f$ contains 
a $\psi_{s}$ for some $s \in ~\C$ where $\psi_s$
is the Bessel function defined in the introduction.
Our aim in this section
is to prove a similar result for the complex semisimple
Lie groups. Our definition of {\it radiality} is taken from
\cite{VV} and it coincides with the definition in the 
previous section
when the function is in $L^1(K \backslash G/K).$

Throughout this section we assume that $G$ is a complex 
semisimple Lie group. Let $Exp:~{\bf p} \rightarrow G/K$ denote
the map $P \rightarrow (exp P)K.$ Then $Exp$ is a 
diffeomorphism. If $dx$ denotes the $G-$invariant measure on
$G/K,$ then

\begin{equation}
\int_{G/K}~f(x)~dx = \int_{\bf p}~f(Exp P)~J(P)~dP,
\end{equation} 
where $$J(P) = det \left (  \frac{\sinh adP}{adP}\right ). $$
Since $G$ is a complex group, the elementary spherical functions
are given by a simple formula:

\begin{equation}
\varphi_\lambda (Exp P) = J(P)^{-\frac{1}{2}}~\int_{K} e^{i
\langle A_\lambda, Ad(k)P \rangle}~dk, \quad P \in {\bf p}.
\end{equation}
Here $A_\lambda$ is the unique element in ${\bf a}_{~\C}$
such that $\lambda(H) = \langle A, A_\lambda \rangle$ for
all $H \in {\bf a}_{~\C}.$

Let $E(K \backslash G/K)$ be the dual of $C^\infty(K \backslash
G /K).$ Then $E(K \backslash G/K)$ can be identified with the
space of compactly supported $K$-biinvariant distributions on
$G/K.$ If $w$ is such a distribution then $\hat{w}(\lambda)
= w(\varphi_\lambda)$ is well defined and is called the 
spherical Fourier transform of $w.$ By the Paley-Wiener theorem
we know that $\lambda \rightarrow \hat{w}(\lambda)$ is an entire
function of exponential type. Similarly, $E(\R^l)$ will denote
the space of compactly supported distribution on $\R^l$ and
$E^W(\R^l)$ consists of the 
Weyl group invariant ones. From the work in \cite{BS79} we know
that the Abel transform $$ S : E(K \backslash G/K) \rightarrow
E^W(\R^l)$$ is an isomorphism and $\widetilde{S(w)}(\lambda) =
\hat{w}(\lambda)$ for $w \in E(K \backslash G/K),$ where $
\widetilde{S(w)}(\lambda)$ is the Euclidean Fourier transform
of the distribution $S(w).$  We also need the following result 
from \cite{BS79}.

\begin{Proposition}
There exists a linear topological isomorphism $T$ from $C^\infty
(K \backslash G/K)$ onto $C^\infty(\R^l)^W$ such that $$S(w)
(T(f)) = w(f)$$ for all $w \in E(K \backslash G/K)$ and $f \in
C^\infty(K \backslash G/K).$ We also have, $$S(w^{'})* T(w \ast
f) = T(w^{'} \ast w \ast f)$$ for all 
$w, w^{'} \in E(K \backslash G/K)$ and
$f \in C^\infty (K \backslash G/K).$ Moreover, $$ T(\varphi_{
\lambda}) = \frac{1}{|W|}~\sum_{\tau \in W} exp (i 
\langle \tau. \lambda, x \rangle).$$

\end{Proposition}

A $K$-biinvariant function  $f$ is called 
{\it radial} if it is of the form $$f(x) = J(Exp^{-1} x)^{
-\frac{1}{2}} u(d (0, x)),$$ where $d$ is the Riemannian 
distance induced by the
the Killing form on $G/K$ and $u$ is a function on 
$[0, \infty).$ Theorem 4.6 in \cite{VV} shows that
this definition of {\it radiality} coincides with the one
in the previous section if the function is integrable.
That is, $f \in L^1(K \backslash G/ K)$ has the above
form if and only if the spherical Fourier transform 
$\hat{f} (\lambda)$
depends only on $(\lambda_1^2 + \lambda_2^2 \cdots + 
\lambda_l^2)^{\frac{1}{2}}.$ We denote the class of smooth
radial functions by $C^\infty(K \backslash G/K)_{rad}$ and
$C_c^\infty(K \backslash G/K)_{rad}$ will consists of 
compactly supported functions in $C^\infty(K \backslash G/K)_{
rad}$.

For $f \in C^\infty(K \backslash G /K)$ define $$f^{\#} 
(Exp P)
= J(P)^{-\frac{1}{2}}~\int_{SO({\bf p})}~J(
\sigma. P)^{\frac{1}{2}}~
f(\sigma. P)~d\sigma,$$ where $SO( {\bf p})$ is the special
orthogonal group on ${\bf p}$ and $d\sigma$ is the Haar measure
on $SO({\bf p}).$ Here, by $f(P)$ we mean $f(Exp P).$ Clearly,
$f \rightarrow f^\#$ is the projection from $C^\infty
(K \backslash G/ K)$ onto $C^\infty(K \backslash G/K)_{rad}.$

\begin{Proposition}
(a) The space $C^\infty (K \backslash G/K)_{rad}$ is reflexive.

\noindent
(b) The strong dual $E(K \backslash G/K)_{rad}$ of $C^\infty(
K \backslash G/K)_{rad}$ is given by $$\{ w \in E(K \backslash
G/K):~\hat{w}(\lambda)~\mbox{is a function of}~
(\lambda_1^2 + \lambda_2^2 + \cdots \lambda_l^2)^{\frac{1}{2}}
\} .$$

\noindent
(c) The space $C^\infty(K \backslash G/K)_{rad}$ is invariant
under convolution by $w \in E(K \backslash G/K)_{rad}$. 
\end{Proposition}

\noindent
{\bf Proof:}~(a) The space $C^\infty (K \backslash G/K)_{rad}$
is a closed subspace of $C^\infty (K \backslash G/K)$ which is
a reflexive Fr$\acute e$chet space.

\noindent
(b) Define $B_\lambda = \varphi_{\lambda}^\#,$ 
the projection 
of $\varphi_\lambda$ into $C^\infty(K \backslash G/K)_{rad}.$ 
A simple computation shows that $$B_\lambda(Exp P) =
J(P)^{-\frac{1}{2}}~\int_{SO({\bf p})}~e^{i \langle A_\lambda,
\sigma. P \rangle}~d\sigma.$$ It is clear that, $B_\lambda$
as a function of $\lambda$ depends only on $(\lambda_1^2 +
\lambda_2^2+ \cdots \lambda_l^2)^{\frac{1}{2}}.$
Now, let $w \in E(K \backslash G/K).$ Define
a distribution $w^\#$ by $w^\# (f) = w( f^\# ).$ 
It is easy to see
that $w^\#$ is a compactly supported $K-$biinvariant 
distribution. Clearly, if $w \in E(K \backslash G/K)_{rad},$
then $w = w^\#.$ It follows that $\hat{w}(\lambda) = w( 
\varphi_\lambda) = w(B_\lambda).$ Consequently, $\hat{w}(
\lambda)$ is a function of $(\lambda_1^2 + \lambda_2^2+ \cdots
+ \lambda_l^2)^{\frac{1}{2}}.$ It also follows that $E(K
\backslash G/K)_{rad}$ is reflexive.

\noindent
(c) Observe that if $w \in E(K \backslash G/K)_{
rad}$ and $g \in C_c^\infty(K \backslash G/K)_{rad}$ then
$w \ast g \in C_c^\infty(K \backslash G/K)_{rad}.$ This follows
from (b) above and Theorem 4.6 in \cite{VV}.
 Next, if $g$ is arbitrary,
we may approximate $g$ with $g_n \in C_c^\infty(K \backslash G/
K)_{rad}.$

We are in a position to state our main result in this section.
Let $V$ be a closed subspace of $C^\infty(K \backslash G/K)_{
rad}.$ We say, $V$ is an ideal in $C^\infty(K \backslash G/K)_{
rad}$ if $f \in V$ and $w \in E(K \backslash G/K)_{rad}$ implies
that $w \ast f \in V.$ 
\bt
(a) If $V$ is a non zero ideal in $ C^\infty(K \backslash G/K)_{
rad}$ then there exists a $\lambda \in {\bf a}_{~\C}^*$
such that $B_\lambda \in V.$

\noindent
(b) If $f \in C^\infty(K \backslash G/K)_{rad},$ then the closed
left $G$ invariant subspace generated by $f$ in $C^\infty( 
G/K)$ contains a $\varphi_\lambda$ for some $\lambda
\in {\bf a}_{~\C}^*.$ 

\et

\noindent
{\bf Proof:}~We closely follow the arguments in \cite{BS79}. 

\noindent
(a) Notice that the map 
$$S:~E(K \backslash G/K)_{rad} \rightarrow E(\R^l)_{rad}$$ 
is a linear topological isomorphism.
Using the 
reflexivity of the spaces involved and arguing as in \cite{BS79}
we obtain that (as in Proposition 3.1)
$$T :C^\infty(K \backslash G/K)_{rad} \rightarrow
C^\infty(\R^l)_{rad}$$ is a linear topological isomorphism, 
where $C^\infty(\R^l)_{rad}$ stands for the space of $C^\infty$
radial functions on $\R^l$ and $$S(w)(T(f) = w(f)~ \forall
w \in E(K \backslash G/K)_{rad}, f \in C^\infty(K \backslash 
G/K)_{rad}.$$

Another application of Proposition 3.1 implies 
that we have a one-one correspondence between the ideals
in $C^\infty(K \backslash G/K)_{rad}$ and $C^\infty(\R^l)_{rad}
.$ Here, ideal in $C^\infty(\R^l)_{rad}$ means a closed subspace
invariant under convolution by compactly supported radial
distributions on $\R^l.$ From \cite{BS90} or \cite{BST}
we know that any ideal
in $C^\infty(\R^l)_{rad}$ contains a $\psi_s$ (Bessel function)
for some $s \in ~\C.$ To complete the proof it suffices to show
that under the topological isomorphism $T$ the function
$B_\lambda$ is mapped into $\psi_s$ where $s^2 = (\lambda_1^2
+ \lambda_2^2 + \cdots + \lambda_l^2)^2.$

Now, we have $S(w)(T(B_\lambda)) = w(B_\lambda).$ Since $w \in
E(K \backslash G/K)_{rad}$ we know that $w(B_\lambda)$
is nothing but $w(\varphi_\lambda)$ which equals $\widetilde{
(Sw)} (\lambda).$ Since $S$ is onto, this implies that $T(
B_\lambda) = \psi_s$ where $s^2 = (\lambda_1^2 + \lambda_2^2
+ \cdots \lambda_l^2)^{\frac{1}{2}}.$

\noindent
(b)  From \cite{BS79} we know that $T(\varphi_\lambda) = 
\psi_\lambda$
where $\psi_\lambda (x) = \frac{1}{|W|}~\sum_{\tau \in W}~
exp (i \tau \lambda .x) .$ Let $V_f$ denote the left 
$G$-invariant
subspace generated by $f.$ Then $T(V_f)$ surely contains the
space $$V_{T(f)} = \{ S(w) \ast T(f):~w \in E(K \backslash G/K) 
\}.$$ From Proposition 3.2,
$T(f)$ is a radial $C^\infty$ function
on $\R^l.$ Hence, from \cite{BST}, the translation invariant 
subspace $X_{T(f)},$ generated by $T(f)$ in $C^\infty(\R^l)$
 contains
a $\psi_s$ for some $s \in ~\C$ and consequently all the
exponentials $e^{i z.x}$ where $z = (z_1, z_2, \dots z_l)$
satisfies $z_1^2 + z_2^2 + \cdots + z_l^2 = s^2.$ Now, it is 
easy to see that the map $g \rightarrow g^W$ where $g^W(x)
= \frac{1}{|W|}~\sum_{\tau \in W}~g(\tau .x),$ from $X_{
T(f)}$ into $V_{T(f)}$ is surjective. Hence, there exists
a $\lambda \in ~\C^l$ such that $\psi_\lambda \in V_{T(f)}.$
Since $T(\varphi_\lambda) = \psi_\lambda,$ this finishes 
the proof.

Our next result is a Wiener-Tauberian type theorem for compactly
supported distributions. Let $E(G/K)$ denote the space of 
compactly supported supported distributions on $G/K.$ If $g \in
G$ and $w \in E(G/K)$ then the left $g-$translate of $w$ is
the compactly supported distribution ${}^{g}w$ defined by 
$${}^g w (f) = w({}^{g^{-1}}f), \quad f \in C^\infty
(G/K)$$ where ${}^{x}f (y) = f(x^{-1}y).$ 

\bt
Let $ \{ w_\alpha:~\alpha \in I \}$ be a family of distributions
contained in $E(K \backslash G/K)_{rad}.$ Then, the left 
$G-$translates 
of this family spans a dense subset of $E(G/K)$ if
and only if there exists no $\lambda \in {\bf a}_{~\C}^*$
such that $\hat{w}_{\alpha} (\lambda) = 0 $ for all $\alpha
\in I.$
\et

\noindent
{\bf Proof:}~We start with the {\it if} part of the
theorem. Let $J$ stand for the 
closed span of the left $G-$translates 
of the distributions $w_\alpha$ in $E(G/K).$ It
suffices to show that $E(K \backslash G/K) \subset J.$ 
To see this, let $f \in C^\infty(G/K)$ be such that $w(f) = 0$
for all $w \in E(K \backslash G/K).$ Since $J$ is left 
$G-$invariant 
we also have $w (f_g) = 0$ for all $g \in G,$ 
where $f_g$ is the $K$-biinvariant function defined by 
$$f_g (x) = \int_K~f(gkx)~dk.$$ It follows that $f_g \equiv
0$ for all $g \in G$ and consequently $f \equiv 0.$

Next, we claim that if $E(K \backslash G/K)_{rad} \subset J$ 
then $E(K \backslash G/K) \subset J.$ To prove this it is enough
to show that $$\{ g \ast w:~w \in E(K \backslash G/K)_{rad}, 
g \in C_c^\infty(K \backslash G/K) \}$$ is dense in
$E(K \backslash G/K).$ Notice that, by Proposition 3.2
the map $S$ from $E(K \backslash G/K)$ onto $E(\R^l)^W$
is a linear topological isomorphism which maps $E(K \backslash
G/K)_{rad}$ onto $E(\R^l)_{rad}$ isomorphically. Hence, it
suffices to prove a similar statement for $E(\R^l)_{rad}$ 
and $E(\R^l)^W$ which is an easy exercise in distribution
theory!

So, to complete the proof of Theorem 3.4 
we only need to show that
$$\{ g \ast w_\alpha:~\alpha \in I, g \in C_c^\infty (K 
\backslash G/K)_{rad} \}$$ is dense in 
$E(K \backslash G/K)_{rad}.$
If not, consider $$J_{rad} = 
\{ f \in C^\infty(K \backslash G/K)_{rad}:
~(g \ast w_\alpha)(f) = 0 ~\forall g \in C_c^\infty(K \backslash
G/K),~\alpha \in I \}.$$ The above is clearly a closed
subspace of $C^\infty(K \backslash G/K)_{rad}$ which is 
invariant under convolution by 
$C_c^\infty(K \backslash G/K)_{rad}.$ 
By Theorem 3.3 we have $B_\lambda \in J_{rad}$ for some $\lambda
\in {\bf a}_{~\C}^*.$ It follows that $\hat{w}_\alpha 
(\lambda) = 0
$ for all $\alpha \in I$ which is a contradiction. This finishes
the proof.

For the {\it only if} part, it suffices to observe that if
$g \in C_c^\infty(G/K)$ then
$$g 
\ast w_\alpha (\varphi_\lambda) = \hat{g^\#}(\lambda) \hat{w}_{
\alpha} (\lambda)$$ where $g^\#(x) = \int_K~g(kx)~dk.$

\noindent
{\bf Remark:}~A similar theorem for {\bf all} rank one spaces 
may be derived from the results in \cite{BS90}.

\vskip.15in
\begin{flushleft}
Department of Mathematics \\
Indian Institute of Science \\
Bangalore -12 \\
India \\
E-mail:~naru@math.iisc.ernet.in, sitaram.alladi@gmail.com \\

\end{flushleft}

\end{document}